\documentclass[12pt]{article}
\usepackage{mathtools} 
\usepackage{amsmath, amssymb} 
\usepackage{amsthm}
\usepackage{enumerate}  
\usepackage{url} 
\usepackage{authblk}
\usepackage{hyperref}

\usepackage{tikz}
\usepackage{tkz-graph}
\usetikzlibrary{shapes}
\usetikzlibrary{arrows}
\usetikzlibrary{decorations.markings}

\usepackage{graphicx}
\usepackage{caption,subcaption}

\newcommand\cx{{\mathbb C}}

\newcommand\ints{{\mathbb Z}}
\newcommand\re{{\mathbb R}}
\newcommand\rats{{\mathbb Q}}

\DeclarePairedDelimiter\abs{\lvert}{\rvert}%
\DeclarePairedDelimiter\norm{\lVert}{\rVert}%

\makeatletter
\let\oldabs\abs
\def\abs{\@ifstar{\oldabs}{\oldabs*}}
\let\oldnorm\norm
\def\norm{\@ifstar{\oldnorm}{\oldnorm*}}
\makeatother

\newcommand\opk[1]{\mathop{\mathrm{#1}}\nolimits}

\newcommand\comp[1]{{\mkern2mu\overline{\mkern-2mu#1}}}
\newcommand\pmat[1]{\begin{pmatrix} #1 \end{pmatrix}}
\newcommand\seq[4]{#1_{#2},#1_{#3},\ldots,#1_{#4}}

\newtheoremstyle{plainsl}%
	{\topsep}
	{\topsep}
	{\slshape} 
	{}
	{\normalfont\bfseries}
	{.}
	{ }
	{}

\swapnumbers

{\theoremstyle{plainsl}
\newtheorem{theorem}{Theorem}[section]
\newtheorem{lemma}[theorem]{Lemma}
\newtheorem{corollary}[theorem]{Corollary}}
{\theoremstyle{remark}
}

\renewcommand\proof{\noindent\textsl{Proof. }}
\newcommand\sqr[2]{{\vbox{\hrule height.#2pt
    \hbox{\vrule width.#2pt height#1pt \kern#1pt
        \vrule width.#2pt}\hrule height.#2pt}}}
\renewcommand\qed{%
	\ifmmode\eqno\sqr53
	\else\nolinebreak\ \hfill\sqr53\medbreak\fi}

\DeclareMathOperator{\col}{col}

\newcommand\ip[2]{\langle#1,#2\rangle}
\newcommand\dist[2]{\opk{dist}(#1,#2)}
\newcommand\one{{\bf1}}

\usepackage{blkarray}

\title{$\epsilon$-Uniform Mixing in Discrete Quantum Walks}
\author{Hanmeng Zhan}
\date{}

\affil{Computer Science Department\\ Worcester Polytechnic Institute, Massachusetts, USA\\\texttt{hzhan@wpi.edu}}

\begin{document}
\maketitle

\begin{abstract}
	We study whether the probability distribution of a discrete quantum walk can get arbitrarily close to uniform, given that the walk starts with a uniform superposition of the outgoing arcs of some vertex. We establish a characterization of this phenomenon on regular non-bipartite graphs in terms of their adjacency eigenvalues and eigenprojections. Using theory from association schemes, we show this phenomenon happens on a strongly regular graph $X$ if and only if $X$ or $\comp{X}$ has parameters $(4m^2, 2m^2\pm m, m^2\pm m, m^2\pm m)$ where $m\ge 2$.

\end{abstract}

\section{Introduction}
Discrete quantum walks are quantum analogues of discrete random walks on graphs, and have been used as subroutines of various quantum algorithms. For example, a quantum search algorithm  \cite{Grover1996,Shenvi2003} is equivalent to a discrete quantum walk, which starts from a uniform superposition of all the arcs of a graph, and gets very close to a state that ``concentrates on" a vertex. 

In this paper, we consider a somewhat related problem: if we start from a state that concentrates on a vertex, can the probability distribution get \textsl{arbitrarily} close to uniform? A similar property, called local $\epsilon$-uniform mixing, has been studied in continuous quantum walks \cite{Godsil2014a, Chan2013, CoutinhoGodsil}, but much less is known in the discrete world. To the best of our knowledge, the instantaneous uniform mixing property of discrete quantum walks has only been analyzed on the hypercubes \cite{Moore2001}. 

To bridge this gap, we develop the general theory for local $\epsilon$-uniform mixing in discrete quantum walks. Using tools from association schemes, we classify all complete graphs and strongly regular graph that admit this phenomenon.

Throughout this paper, all graphs are connected and regular. We will consider arc-reversal quantum walks on these graphs with Grover coins \cite{Kendon2003}. Given a graph $X$, an \textsl{arc} in $X$ is an ordered pair $(a,b)$ of adjacent vertices, with \textsl{tail} $a$ and \textsl{head} $b$. Let $D_t$ denote the arc-tail incidence matrix, and $D_h$ the arc-head incidence matrix. If $X$ is $k$-regular, the transition matrix of our quantum walk is equal to
\[U = R\left(\frac{2}{k} D_t^TD_t - I\right).\]
A state associated with $X$ is a one-dimensional subspace of $\cx^{\mathrm{arcs(X)}}$, represented by a unit vector in it. Given that the walk starts with state $x$, its state at time $t$ is given by $U^tx$. The Schur product $(U^tx)\circ \comp{(U^t x)}$ represents the probability distribution over all arcs of the graph at time $t$.

\section{Spectral decomposition}
We can translate any phenomenon of a quantum walk into properties of the spectral decomposition of $U$; the question here is whether they indicate something special about the underlying graph. Indeed, if the graph is regular, then there is a spectral correspondence between the adjacency matrix $A$ of $X$ and the transition matrix $U$, as shown below.

\begin{theorem}\cite[Ch 3]{Godsil2023} \label{spd}
	Let $X$ be  $k$-regular graph with spectral decomposition
	\[A = \sum_{\lambda} E_{\lambda}.\]
	The the transition matrix $U$ has spectral decomposition
	\[U = 1 \cdot F_1 + (-1) \cdot F_{-1} + \sum_{\theta\in(0,\pi)} (e^{i\theta} F_{\theta} + e^{-i\theta} F_{-\theta}),\]
	where the eigenvalues and eigenspaces are given as follows.
	\begin{enumerate}[(i)]
		\item For the $1$-eigenspace, we have
		\[\col(F_1) = (\col(I+R)\cap\col(D_t^T)) \oplus (\ker(I+R)\cap\ker(D_t)).\]
		Moreover, the first subspace is given by
		\[ (\col(I+R)\cap\col(D_t^T)) = D_t^T \col(E_k),\]
		while the second subspace is given by 
		\[\ker(I+R)\cap \ker(D_t) = N \ker(C),\]
		where $C$ is the signed vertex-edge incidence matrix of any orientation of $X$, and $N$ is the corresponding signed arc-edge incidence matrix.
		\item For the $(-1)$-eigenspace, we have
		\[\col(F_{-1}) = (\ker(I+R)\cap \col(D_t^T))\oplus (\col(I+R)\cap\ker(D_t)),\]
		Moreover, the first subspace is given by 
		\[(\ker(I+R)\cap \col(D_t^T))=D_t^T \ker(E_{-k}),\]
		while the second subspace is given by
		\[(\col(I+R)\cap\ker(D_t))=M\ker(B),\]
		where $B$ is the unsigned vertex-edge incidence matrix, and $M$ is the unsigned arc-edge incidence matrix.
		\item For the $e^{\pm i\theta}$-eigenspace, we have $\lambda=\cos\theta$, and
			\[F_{\pm \theta} = \frac{1}{2k\sin^2\theta} (D_t - e^{\pm i\theta} D_h)^T E_{\lambda} (D_t - e^{\mp i \theta} D_h).\]
	\end{enumerate}

\end{theorem}

Unless otherwise specified, given a graph $X$, we let $A$ denote the adjacency matrix, $D$ the degree matrix, $B$ the unsigned vertex-edge incidence matrix, $M$ the unsigned arc-edge incidence matrix, $C$ the signed vertex-incidene matrix of some orientation of $X$, and $N$ the corresponding signed arc-edge incidence matrix. We list some properties of these matrices that we will use frequently.

\begin{lemma}\label{props}
	The following hold for any graph.
	\begin{enumerate}[(1)]
		\item $D_tR = D_h$ and $D_h R = D_t$.
		\item $RM=M$ and $RN=-N$.
		\item $MM^T = I+R$ and $NN^T=I-R$.
		\item $D_tD_t^T =D_hD_h^T= D$.
		\item $D_t D_h^T = D_hD_t^T = A$.
		\item $D_tM =D_hM =  B$.
		\item $D_tN = -D_h N = C$.
		\item $\col(I+R)\cap \ker(D_t) = M\ker(B)$.
		\item $\ker(I+R)\cap \ker(D_t) = N\ker(C)$.
		\item $(\col(I+R)\cap \ker(D_t))^{\perp} = \col(I-R)+\col(D_t^T)$
		\item $(\ker(I+R)\cap \ker(D_t))^{\perp} = \col(I+R)+\col(D_t^T)$
		\item $\col(I+R) + \col(D_t^T)$ and $\col(I-R) + \col(D_t^T)$ are $R$-invariant.
	\end{enumerate}
\end{lemma}
\proof
Properties (1) and (2) hold as $R$ swaps the tail of each arc with its head. Properties (3) and (4) are immediate from the definition. For (5), let $u$ and $v$ be any two vertices. Then 
\[(D_tD_h^T)_{u,v} = 
\begin{cases}
	1,\quad \text{if $u\sim v$ in X},\\
	0,\quad \text{otherwise}.
\end{cases}\]
For (6) and (7), let $u$ and $e$ be a vertex and an edge. Then
\[(D_tN)_{u, e} = \begin{cases}
	1,\quad \text{if $e$ leaves $u$ in the given orientation of $X$},\\
	-1,\quad \text{if $e$ enters $u$ in the given orientation of $X$},\\
	0,\quad \text{otherwise},
\end{cases}\]
and
\[(D_tM)_{ue} = \begin{cases}
	1,\quad \text{if $u$ is incidence to $e$ in X},\\
	0,\quad \text{otherwise}.
\end{cases}\]
Properties(8) and (9) are from Theorem \ref{spd} but can also be verified using properties (3), (6) and (7). For (10) and (11), simply note that 
\[\col(I+R) = \ker(I-R).\]
To see (12), write $RD_t^T$ as
\[RD_t^T  = (I+R) D_t^T - D_t^T = (I-R)(-D_t^T) + D_t^T \tag*{\sqr53}\]

We will start our quantum walk with a ``uniform superposition of the outgoing arcs of vertex $a$", that is, the state
\[x_a:=\frac{1}{\sqrt{\deg(a)}} D_t^T e_a.\]
We say a matrix is \textsl{real} if all its entries are real, and \textsl{flat} if all its entries have the same absolute value. The following theorem describes the entries of the state after time $t$. 

\begin{theorem}\label{bip-non-bip}
	Let $X$ be a $k$-regular graph. For any time $t$, the $(u,v)$-entry of $U^t x_a$ is
	\begin{equation}\label{eqn:amp}
		\frac{1}{\sqrt{k}}\left(\sum_{\theta_r \in (0,\pi)} \frac{1}{\sin\theta_r}(\sin(t\theta_r) (E_r)_{va} - \sin((t-1)\theta_r) (E_r)_{ua}) + (E_k)_{ua} + (-1)^t(E_{-k})_{ua}\right)
	\end{equation}
	In particular:
	\begin{enumerate}[(i)]
		\item If $X$ is non-bipartite, then at any time $t$, the state $U^t x_a$ is real.
		\item If $X$ is bipartite, then at any time $t$, the state $U^tx_a$ is has flat imaginary part. More precisely, the imaginary part of the $(u,v)$-entry of $U^tx_a$ is 
	\[	\begin{cases}
			\frac{1}{n\sqrt{k}}\sin (t\pi),\quad &\text{ if $u$ and $a$ lie in the same color class}, \\
			-\frac{1}{n\sqrt{k}}\sin (t\pi),\quad & \text{ if $u$ and $a$ lie in different color classes}.
		\end{cases}\]
	\end{enumerate} 
\end{theorem}
\proof
The formula follows from Theorem \ref{spd} and Lemma \ref{props} (4) and (5). Since the only term that could contribute an imaginary number is $\frac{1}{\sqrt{k}}(-1)^k E_{-k}$, the state $U^tx_a$ is real except when $E_{-k}\ne 0$, that is, when $X$ is bipartite. Moreover, for a $k$-regular graph, 
\[E_{-k}=\frac{1}{n}\pmat{J & -J\\-J&J},\]
and so $(E_{-k})_{ua}$ depends only on whether $u$ and $a$ lie in the same color class. Expanding $(-1)^t=\cos(t\pi) + i\sin(t\pi)$ gives the imaginary part of $U^tx_a$. 
\qed

One consequence of Theorem \ref{bip-non-bip} is that for regular non-bipartite graphs, the closure of $\{U^tx_ax_a^TU^{-t}\}_{t\in \re}$ contains only real matrices.

\section{Uniform and $\epsilon$-uniform mixing}
A graph $X$ with $m$ arcs admits \textsl{local uniform mixing} at vertex $a$  if there is a time $t\in\re$ such that
\[(U^t x_a) \circ \comp{(U^t x_a)} = \frac{1}{m}\one.\]
We say $X$ admits \textsl{uniform mixing} if it admits local uniform mixing at every vertex. Below we determine all $K_n$ and $K_{k,k}$ that admit local uniform mixing.

\begin{lemma}
	$K_2$ admits uniform mixing at time $1/2$. $K_4$ admits uniform mixing at time $\pi/\arccos(-1/3)$. No other complete graphs admit local uniform mixing. 
\end{lemma}
\proof
The spectral decomposition of $A(K_n)$ is
\[A(K_n) = (n-1) \cdot \frac{1}{n}J + (-1) \cdot \left(I-\frac{1}{n}J \right).\]
For $n=2$, Equation \eqref{eqn:amp} becomes
\[\frac{1}{2}\pmat{1&1\\1&1}_{ua} + \frac{(-1)^t}{2}\pmat{1 & -1\\-1&1}_{ua}=
\begin{cases}
	\frac{1}{2}(1+(-1)^t),\quad u=a\\
	\frac{1}{2}(1-(-1)^t),\quad u\ne a
\end{cases}.\]
Thus at time $t=1/2$, the $(u,v)$ entry of $U^t x_a$ is 
\[\begin{cases}
	\frac{1}{2}(1+i),\quad \text{ if }u=a\\
	\frac{1}{2}(1-i),\quad \text{ if } u\ne a.
\end{cases}\]
That is, uniform mixing occurs on $K_2$ at time $1/2$.

For $n\ge 3$, Equation \eqref{eqn:amp} becomes
\begin{equation}\label{uv-entry}
	\frac{1}{\sqrt{n-1}}\left(\frac{\sin(t\theta)}{\sin\theta} \left(I-\frac{1}{n}J\right)_{va} - \frac{\sin((t-1)\theta)}{\sin\theta} \left(I - \frac{1}{n}J\right)_{ua} + \frac{1}{n} \right), 
\end{equation}
where $\theta=\arccos(-1/n)$. If $n=4$ and $t=\pi/\arccos(-1/3)$, then Equation \eqref{uv-entry} reduces to
\begin{align*}
	\frac{1}{\sqrt{3}}\left(\frac{1}{4} - \left(I - \frac{1}{4}J\right)_{ua} \right) = \begin{cases}
		-\frac{1}{\sqrt{12}},\quad &\text{ if }u=a\\
		\frac{1}{\sqrt{12}},\quad & \text{ if }u\ne a,
	\end{cases}
\end{align*}
which shows uniform mixing on $K_4$. For other values of $n$, let $a, b, c$ be three distinct vertices in $K_n$, and suppose local uniform mixing occurs at vertex $a$. If either of $e_{(a,b)}^TU^t x_a$ and $e_{(b,a)}^TU^t x_a$ equals $e_{(b,c)}^TU^t x_a$, then by Equation \eqref{uv-entry}, one of $\sin(t\theta)$ and $\sin((t-1)\theta)$ is zero, and so 
\[\frac{2}{n}=\frac{1}{\sqrt{n}}.\]
Hence it must be that 
\[e_{(a,b)}^TU^t x_a= e_{(b,a)}^TU^t x_a=- e_{(b,c)}^TU^t x_a=\pm \frac{1}{\sqrt{n(n-1)}}.\]
However, this has no solution according to the formula in Equation \eqref{uv-entry}.
\qed

\begin{lemma}\label{lem: nonbip_um}
	$K_{k,k}$ admits local uniform mixing if and only if $k=1$. 
\end{lemma}
\proof
Consider $K_{k,k}$ where $k\ge 2$. The spectral decomposition of $A(K_{k,k})$ is
\[A(K_{k,k}) = k\cdot \frac{1}{2k}\pmat{J & J \\J& J} + (-k)\cdot \frac{1}{2k} \pmat{J & -J\\ -J & J} + 0\cdot \pmat{I-\frac{1}{k} J & 0 \\ 0 & I-\frac{1}{k}J}\]
By Equation \eqref{eqn:amp}, the $(u,v)$-entry of $U^tD_t^T e_a$ is
\begin{align*}
	e_{(u,v)}^TU^t D_t^T e_a =&\sin(t\theta) \pmat{I-\frac{1}{k} J & 0 \\ 0 & I-\frac{1}{k}J}_{va} +\cos(t\theta) \pmat{I-\frac{1}{k} J & 0 \\ 0 & I-\frac{1}{k}J}_{ua} \\
	&+\frac{1}{2k}\pmat{J & J \\J& J} _{ua} + (-1)^t \frac{1}{2k} \pmat{J & -J\\ -J & J}_{ua}
\end{align*}
where $\theta=\pi/2$. In particular, if $u$ lies in the same color class as $a$, then 
\begin{equation}\label{eqn: cases}
	e_{(u,v)}^TU^t D_t^T e_a=
	\begin{cases}
		\frac{1}{k} \left(\cos^2(t\theta)+(k-1)\cos(t\theta) + i\sin(t\theta)\cos(t\theta)\right),\quad \text{ if } u= a\\
		\frac{1}{k} \left(\cos^2(t\theta)-\cos(t\theta) + i\sin(t\theta)\cos(t\theta)\right),\quad \text{ if } u\ne a\\
	\end{cases}
\end{equation}
Suppose local uniform mixing occurs at vertex $a$ at time $t$. Then for any arc $(u,v)$ we have
\[\norm{e_{(u,v)}^TU^t D_t^T e_a }^2 = \frac{1}{2k}\]
It follows from Equation \eqref{eqn: cases} that $\cos(t\theta)\ne 0$, and the two cases have opposite real parts. Solving for $\cos(t\theta)$ and then $k$ yields $(k-2)^2 = 2$, which is impossible.
\qed

We turn our attention to relaxations of uniform mixing and local uniform mixing. A graph $X$ with $m$ arcs admits \textsl{local $\epsilon$-uniform mixing} at vertex $a$ if for any $\epsilon>0$, there is a time $t\in\re$ such that
\[\norm{(U^t x_a) \circ \comp{(U^t x_a)} - \frac{1}{m}\one } < \epsilon.\]
It admits \textsl{$\epsilon$-uniform mixing} if local $\epsilon$-uniform mixing occurs at every vertex. 

Our goal is to determine all complete graphs and strongly regular graphs that admit local $\epsilon$-uniform mixing. We settle the bipartite cases in this section.

\begin{lemma}\label{bip_srg}
	For $k\ge 2$, the complete bipartite graph $K_{k,k}$ does not admit local $\epsilon$-uniform mixing.
\end{lemma}
\proof
Since the eigenvalues of $A(K_{k,k})$ are $k$, $-k$ and $0$, the eigenvalues of $U$ are $\pm 1, \pm i$. Therefore  $U^4=I$, and for any $t\in \re$,
\[(U^{t+4}x_a) \circ \comp{(U^{t+4} x_a)} = (U^t x_a) \circ \comp{(U^t x_a)}. \]
Thus, if $(U^t x_a) \circ \comp{(U^t x_a)} $ gets arbitrarily close to the uniform probability distribution, then it must attain this distribution at some $t\in[0,4]$. However, Lemma \ref{lem: nonbip_um} shows that $K_{k,k}$ does not admit local uniform mixing for $k\ge 2$.
\qed

\section{Perfect and pretty good state transfer}
One important question in quantum walks is to decide whether the system can get arbitrary close to a target state. We say a graph admits perfect state transfer from $x$ to $y$ if there is a time $t\in \re$ and a unimodular $\gamma\in \cx$ such that 
\[U^t x = \gamma y.\]
We say there is \textsl{pretty good state transfer} from state $x$ to state $y$ if there exists a unimodular $\gamma \in \cx$ such that, for any $\epsilon >0$, there is a time $t\in\re$ with 
\[\norm{U^tx - \gamma y}<\epsilon.\]
Clearly, local uniform mixing is equivalent to perfect state transfer from $x_a$ to a flat state. In the following lemma, we show that local $\epsilon$-uniform mixing is also a special case of pretty good state transfer.

\begin{lemma}
	A graph with $m$ arcs admits local $\epsilon$-uniform mixing at vertex $a$ if and only if there is pretty good state transfer from $x_a$ to some flat state $y\in \cx^m$. 
\end{lemma}
\proof
The backward direction is clear: as each entry of $U^tx_a$ is bounded,  there is a constant $c$ such that
\[\norm{(U^t x_a) \circ \comp{(U^t x_a)} - y\circ \comp{y}} \le  c \norm{U^t x_a - \gamma y}.\]
To prove the forward direction, assume $X$ admits local $\epsilon$-uniform mixing at $a$. Then for any $\epsilon>0$, there is a $t\in \re$ satisfying
\[\norm{(U^t x_a) \circ \comp{(U^t x_a)} - \frac{1}{m}\one } < \epsilon.\]
Hence there is a sequence of times $\{t_1, t_2, t_3, \cdots\}$ such that 
\[\lim_{j\to \infty} (U^{t_j} x_a) \circ \comp{(U^{t_j} x_a)} =  \frac{1}{m}\one. \]
It follows that $\{U^{t_j} x_a\}_j$ is bounded and has a convergent subsequence $\{U^{t_{j_{\ell}}} x_a \}_{\ell}$, say with limit $y$. As 
\[y\circ \comp{y} = \lim_{\ell\to \infty} (U^{t_{j_{\ell}}}x_a)\circ \comp{(U^{t_{j_{\ell}}}x_a)} =\lim_{j\to \infty} (U^{t_j} x_a) \circ \comp{(U^{t_j} x_a)} =  \frac{1}{m}\one,\]
the entries of $y$ must have absolute value $1/\sqrt{m}$.
\qed

We can say more if the graph is regular, using results from Theorem \ref{bip-non-bip}.

\begin{corollary}
	A regular graph admits local $\epsilon$-uniform mixing at vertex $a$ if and only if there is pretty good state transfer from $x_a$ to some state $y$ where both $\mathrm{Re}(y)$ and $\mathrm{Im}(y)$ are flat. Moreover, if this graph is non-bipartite, then $y$ is a scalar multiple of a $\{\pm 1\}$-vector.
\end{corollary}

Pretty good state transfer on discrete quantum walks has been studied in \cite{Chan2023} in a general setting. In particular, there is a characterization using the eigenvalues and eigenprojections of the transition matrix. We call two states $x$ and $y$ \textsl{strongly cospectral} with respect to $U$ if for each eigenprojection $F_{\theta}$ of $U$, there is $\delta\in\re$ such that
\[F_{\theta} x = e^{i\delta} F_{\theta} y.\]

\begin{theorem}\cite{Chan2023} \label{char-pgst}
Let $U$ be the transition matrix of a quantum walk. There is pretty good state transfer from $x$ to $y$ if and only if the following hold.
	\begin{enumerate}[(i)]
		\item $x$ and $y$ are strongly cospectral: for each eigenprojection $F_{\theta}$ of $U$, there is some $\delta\in\re$ such that $F_{\theta} x = e^{i\delta} F_{\theta} y$.
		\item There is a unimodular $\gamma\in \cx$ such that, for any $\epsilon>0$, there exists $t\in \re$ so that when $F_{\theta} x\ne 0$,
		\[\abs{e^{i(t\theta+\delta)} - \gamma}<\epsilon.\]
	\end{enumerate}
\end{theorem}

We will study Conditions (i) and (ii) separately, in Sections \ref{sec:sc} and \ref{sec: Ka}.

\section{Strong cospectrality \label{sec:sc}}

In this section, we characterize real states that are strongly cospectral to $x_a$ for regular non-bipartite graphs. We start by translating Theorem \ref{char-pgst} (i) into conditions on the spectral decompositon of the adjacency matrix.

\begin{lemma}
	Let $X$ be a $k$-regular non-bipartite graph on $n$ vertices. Let
	\[A = \sum_r k\cos\theta_r E_r\]
	be the spectral decomposition of the adjacency matrix of $X$, assuming $\theta_0=0$. Let $a$ be a vertex in $X$. Let 
	\[U = 1 \cdot F_1 + (-1) \cdot F_{-1} + \sum_{r\ne 0} (e^{i\theta_r} F_{\theta_r} + e^{-i\theta_r} F_{-\theta_r})\]
	be the spectral decomposition of the transition matrix. Then for any $y\in \re^{kn}$, the following hold.
	\begin{enumerate}[(i)]
		\item $F_1 x_a = e^{i\alpha} F_1 y$ for some $\alpha \in \re$ if and only if 
		\[y\in \col(I+R) + \col(D_t^T) \quad \mathrm{and} \quad \sqrt{k} E_0 e_a = \pm E_0 D_t y=\pm E_0 D_h y.\]
		\item $F_{-1} x_a = e^{i\beta} F_{-1} y$ for some $\beta\in \re$ if and only if 
		\[y\in \col(I-R) + \col(D_t^T)\]
		\item For each $r\ne 0$, the condition $F_{\theta_r} x_a = e^{i\delta_r} F_{\theta} y$ is equivalent to $F_{-\theta_r} x_a = e^{-i\delta_r} F_{-\theta_r} y$, and it holds for some $\delta\in \re$ if and only if 
		\[
			\sqrt{k}\cos \delta_r E_r e_a  = E_r D_t y \quad \mathrm{and} \quad 
			\sqrt{k}\cos (\delta_r + \theta_r)E_r e_a  = E_r D_h y
		\]
	\end{enumerate}
\end{lemma}
\proof
We prove each part of the statement.
\begin{enumerate}[(i)]
	\item By Theorem \ref{spd}, $F_1$ is real, so $F_1 x_a = e^{i\alpha} F_1 y$ if and only if $F_1 x_a = \pm F_1 y$, that is,
	\[F_1 (x_a \mp y) =0.\]
	Since the $1$-eigenspace of $U$ is a direct sum
	\[(D_t^T \col(E_0))\oplus (\ker(I+R)\cap\ker(D_t)),\]
	$F_1 (x_a \mp y) =0$ if and only if each of $D_t^T \col(E_0)$ and $\ker(I+R)\cap\ker(D_t)$ is orthogonal to $x_a \mp y$. The first orthogonality condition is equivalent to
	\[\sqrt{k} E_0e_a = \pm E_0 D_t y =\pm E_0 D_h y,\]
	and the second orthogonality condition is equivalent to 
	\[y\in (\ker(I+R)\cap\ker(D_t))^{\perp} = \col(I+R) + \col(D_t^T).\]
	\item Since $X$ is non-bipartite, by Theorem \ref{spd}, the $(-1)$-eigenspace of $U$ is 
	\[\col(I+R)\cap \ker(D_t),\]
	which is orthogonal to $x_a$. Hence $F_{-1} x_a = e^{i\beta} F_{-1}y$ if and only if 
	\[F_{-1}y=0,\]
	that is, 
	\[y\in (\col(I+R)\cap \ker(D_t))^{\perp} = \col(I-R) + \col(D_t^T).\]
	\item The first claim follows from the fact that $F_{-\theta_r} = \comp{F_{\theta_r}}$. For the second claim, we compute $F_{\theta_r} x_a$ and $F_{\theta_r}y$ using the formula in Theorem \ref{spd}:
	\begin{align*}
		F_{\theta_r} x_a &=\frac{i e^{-i\theta_r}}{2\sqrt{k}\sin\theta_r} (D_t - e^{i\theta_r} D_h)^T E_r e_a\\
		F_{\theta_r} y &= \frac{1}{2k\sin^2\theta_r} (D_t-e^{i\theta_r} D_h)^T E_r (D_t - e^{-i\theta_r} D_h)y.
	\end{align*}
Note that
\[D_t-e^{i\theta_r} D_h = D_t(I-e^{i\theta_r} R),\]
which has full row rank as $e^{i\theta_r}\ne \pm 1$. Thus $F_{\theta_r} x_a = e^{i\delta_r} F_{\theta_r} y$ if and only if
\begin{equation}\label{eqn: vxsc}
	i\sqrt{k} \sin \theta_r E_r e_a = e^{i(\delta_r+\theta_r)} E_r D_ty - e^{i\delta_r} D_h y.
\end{equation}
Comparing the real parts of both sides shows
\[\cos(\delta_r +\theta_r) E_r D_t y =\cos \delta_r E_r D_h y.\]
Using Equation \eqref{eqn: vxsc} again, we obtain
\[
\sqrt{k}\cos \delta_r E_r e_a  = E_r D_t y \quad \mathrm{and} \quad 
\sqrt{k}\cos (\delta_r + \theta_r)E_r e_a  = E_r D_h y.
\tag*{\sqr53}\]
	
\end{enumerate}

Putting all these together yields the following characterization of real states strongly cospectral to $x_a$.

\begin{theorem}\label{strcosp}
	Let $X$ be a $k$-regular non-bipartite graph on $n$ vertices. Let
	\[A = \sum_r k\cos\theta_r E_r\]
	be the spectral decomposition of the adjacency matrix of $X$, assuming $\theta_0=0$. Let $a$ be a vertex in $X$. Then $x_a$ is strongly cospectral to a real state $y$ if and only if 
	\begin{enumerate}[(i)]
		\item $y$ lies in the intersection of $\col(I+R)+\col(D_t^T)$ and $\col(I-R) + \col(D_t^T)$, and
		\item there exist real numbers $\delta_r$ where $\delta_0\in \{0, \pi\}$ such  that
			\begin{align*}
			D_t y &=\sqrt{k}\sum_{r} \cos\delta_r E_re_a,  \\
			D_h y &= \sqrt{k}\sum_{r} \cos(\delta_r + \theta_r) E_re_a.
		\end{align*}
	\end{enumerate}
\end{theorem}

We immediately obtain a number-theoretic constraint on the size of the graph when $x_a$ is strongly cospectral to some real flat state.

\begin{corollary}\label{perfect_square}
	Let $X$ be a regular non-bipartite graph on $n$ vertices. If for some vertex $a$, the state $x_a$ is strongly cospectral to a real flat state, then $n$ is a perfect square.
\end{corollary}
\proof
By Theorem \ref{strcosp} (ii),
\[E_0 D_t y =\pm \sqrt{k} E_0 e_a,\]
which implies 
\[\ip{\one}{y}=\pm \sqrt{k}.\]
Since $y$ is flat, its entries are $\pm 1/\sqrt{kn}$. Plugging this into the equation above shows that $\sqrt{n}$ is an integer.
\qed

Another consequence that will become helpful later is that real states strongly cospectral to $x_a$ come in pairs.

\begin{corollary}\label{sc_R}
	Let $X$ be a $k$-regular non-bipartite graph. If $x_a$ is strongly cospectral to a real state $y$, then $x_a$ is strongly cospectral to $Ry$.
\end{corollary}
\proof
First, since $\col(I+R)+\col(D_t^T)$ and $\col(I-R)+\col(D_t^T)$ are $R$-invariant, it contains $y$ if and only if it contains $Ry$. Next, note that $D_tR = D_h$, and so with $\delta_r'=-\delta_r-\theta_r$, the identities
\begin{align*}
	D_t y &=\sqrt{k}\sum_{r} \cos\delta_r E_re_a,  \\
	D_h y &= \sqrt{k}\sum_{r} \cos(\delta_r + \theta_r) E_re_a
\end{align*}
are equivalent to 
			\begin{align*}
	D_t Ry &=\sqrt{k}\sum_{r} \cos\delta_r' E_re_a,  \\
	D_h Ry &= \sqrt{k}\sum_{r} \cos(\delta_r' + \theta_r) E_re_a. 
\end{align*}
Therefore by Theorem \ref{strcosp}, $x_a$ is strongly cospectral to $y$ if and only if it is strongly cospectral to $Ry$.
\qed

At the end of this section, we consider a special form of strong cospectrality, between $x_a$ and real states in $\col(D_t^T)$. By the above result, this also gives a characterization of strong cospectrality between $x_a$ and real states in $\col(D_h^T)$.

\begin{corollary}\label{in-Dt}
Let $X$ be a $k$-regular non-bipartite graph on $n$ vertices. Let
\[A = \sum_r k\cos\theta_r E_r\]
be the spectral decomposition of the adjacency matrix of $X$, assuming $\theta_0=0$. Let $a$ be a vertex in $X$. Then $x_a$ is strongly cospectral to $D_t^Tw$ for some $w\in \re^n$ if and only if 
\[
	w =\frac{1}{\sqrt{k}}\sum_{r} (-1)^{\sigma_r} E_r e_a
\]
for some integers $\{\sigma_r\}$.
\end{corollary}
\proof
We substitute $y=D_t^Tw$ into Theorem \ref{strcosp} (i) and (ii). Condition (i) is vacuously satisfied. Condition (ii) is equivalent to
\begin{align}
	kw &= \sqrt{k} \sum_{r} \cos\delta_r E_r e_a \label{kw}\\
	Aw &=\sqrt{k} \sum_{r} \cos(\delta_r+\theta_r)  e_a.\label{Aw}
\end{align}
Given
\[
w =\frac{1}{\sqrt{k}}\sum_{r} (-1)^{\sigma_r} E_r e_a,
\]
Equations \eqref{kw} and \eqref{Aw} hold with $\delta_r = \pi \sigma_r$. Conversely, given Equations \eqref{kw} and \eqref{Aw}, we have
\[\sum_{r} \cos\delta_r \lambda_r E_r e_a = \sqrt{k}Aw = \sum_{r} \cos(\delta_r + \theta_r) k E_re_a.\]
Hence for each $r$, 
\[\cos\delta_r \cos\theta_r E_r e_a = \cos(\delta_r + \theta_r) E_r e_a.\]
If $E_r e_a\ne 0$, then 
\[\sin\delta_r \sin \theta_r =0,\]
and since $\cos\theta_r\ne \pm1$, we must have $\sin \delta_r=0$. If $E_r e_a =0$, then we may choose any value for $\delta_r$, say $\delta_r =0$. In both cases, $\cos\delta_r =\pm1=(-1)^{\sigma_r}$ for some integer $\sigma_r$, and so 
\[
	w =\frac{1}{\sqrt{k}}\sum_{r} (-1)^{\sigma_r} E_r e_a. \tag*{\sqr53}
\]

\section{Kronecker's approximation \label{sec: Ka}}
To translate Theorem \ref{char-pgst} (ii) into a condition on the adjacency spectrum, we cite the two versions of Kronecker's approximation theorem. The first one helps us characterize pretty good state transfer.

\begin{theorem}\cite[Ch 3]{Levitan1983}
	Given $\seq{\alpha}{1}{2}{n}\in \re$ and $\seq{\beta}{1}{2}{n}\in \re$, the following are equivalent.
	\begin{enumerate}[(i)]
		\item For any $\epsilon>0$, the system
		\[\abs{q\alpha_r - \beta_r}<\epsilon \pmod{2\pi}, \quad r=1,2,\cdots, n\]
		has a solution $q\in\re$.
		\item For any set $\{\seq{\ell}{1}{2}{n}\}$ of integers such that 
		\[\ell_1 \alpha_1 + \cdots + \ell_n \alpha_n =0,\]
		we have
		\[\ell_1 \beta_1 + \cdots + \ell_n \beta_n \equiv 0\pmod{2\pi}.\]
	\end{enumerate}
\end{theorem}

The second version helps us determine a streghthening of pretty good state transfer, where we can approximate the target state in integer times. 

\begin{theorem}\cite{Gonek2016}\label{thm:Kronecker}
	Given $\seq{\alpha}{1}{2}{n}\in \re$ and $\seq{\beta}{1}{2}{n}\in \re$, the following are equivalent.
	\begin{enumerate}[(i)]
		\item For any $\epsilon>0$, the system
		\[\abs{q\alpha_r - \beta_r-p_r}<\epsilon, \quad r=1,2,\cdots, n\]
		has a solution $\{q, \seq{p}{1}{2}{n}\} \in \ints^{n+1}$.
		\item For any set $\{\seq{\ell}{1}{2}{n}\}$ of integers such that 
		\[\ell_1 \alpha_1 + \cdots + \ell_n \alpha_n \in \ints,\]
		we have
		\[\ell_1 \beta_1 + \cdots + \ell_n \beta_n \in \ints.\]
	\end{enumerate}
\end{theorem}

We now give a characterization of pretty good state transfer from $x_a$ to $y$ on regular non-bipartite graphs.

\begin{theorem}\label{thm:pgst}
	Let $X$ be a $k$-regular non-bipartite graph on $n$ vertices. Let
	\[A = \sum_r k\cos\theta_r E_r\]
	be the spectral decomposition of the adjacency matrix of $X$, assuming $\theta_0=0$. Let $a$ be a vertex in $X$. Then $X$ admits pretty good state transfer from $x_a$ to $y$ if and only if $y$ is a scalar multiple of a real state and all of the following hold.
		\begin{enumerate}[(i)]
		\item $y$ lies in the intersection of $\col(I+R)+\col(D_t^T)$ and $\col(I-R) + \col(D_t^T)$.
		\item There exist real numbers $\delta_r$ where $\delta_0\in \{0, \pi\}$ such  that
		\begin{align*}
			D_t y &=\sqrt{k} \sum_{r} \cos\delta_r E_re_a,  \\
			D_h y &= \sqrt{k}  \sum_{r} \cos(\delta_r + \theta_r) E_re_a
		\end{align*}
		\item For any set $\{\ell_r: r\in \Phi_a\}$ of integers such that
		\begin{equation}\label{weak_pgst}
			\sum_{r\in \Phi_a} \ell_r \theta_r =0,
		\end{equation}
		we have
		\[\sum_{r\in \Phi_a} \ell_r\delta_r \equiv 0 \pmod{2\pi}.\]
	\end{enumerate}
	Moreover, if (i) - (iii) hold with Equation \eqref{weak_pgst} replaced by
	\[\sum_{r\in \Phi_a} \ell_r \theta_r \equiv 0 \pmod{2\pi},\]
	then $y$ can be approximated in integer steps.
\end{theorem}
\proof
By scaling $y$ if necessary, we may assume that $y$ is real and $\ip{\one}{y}\ge 0$. Thus Conditions (i) and (ii) are equivalent to strong cospectrality between $x_a$ and $y$. To see why Condition (iii) is equivalent to Theorem \ref{char-pgst} (ii), we note that $e^{it\theta_0} =1$, and $F_{\theta_r} x_a = e^{i\delta_r} F_{\theta_r} y$ if and only if $F_{-\theta_r} x_a = e^{-i\delta_r} F_{-\theta_r} y$ . Hence Theorem \ref{char-pgst} (ii) holds if and only if for any $\epsilon>0$, there exists $t\in \re$ such that 
\[\abs{e^{i(t\theta_r + \delta_r)} -1} <\epsilon,\quad r\in \Phi_a;\]
that is, for any $\epsilon>0$, the system
\[\abs{t\theta_r + \delta_r}<\epsilon \pmod{2\pi}, \quad r\in \Phi_a\]
has a solution $t\in \re$. The result then follows from Kronecker's approximation theorems.
\qed

Using Corollary \ref{sc_R} and its proof, we obtain the following.

\begin{corollary}
	Let $X$ be a regular non-bipartite graph on $n$ vertices. Let $a$ be a vertex in $X$. Then $X$ admits pretty good state transfer from $x_a$ to $y$ if and only if $X$ admits pretty good state transfer from $x_a$ to $Ry$.
\end{corollary}
\proof
Let $\theta_r$, $\delta_r$ and $\ell_r$ be defined as in Theorem \ref{thm:pgst}. If 
\[\sum_{r\in \Phi_a} \ell_r \theta_r \equiv 0\pmod{2\pi},\quad \sum_{r\in \Phi_a} \ell_r\delta_r \equiv 0 \pmod{2\pi},\]
then
\[\sum_{r\in \Phi_a} \ell_r(-\delta_r-\theta_r) \equiv 0 \pmod{2\pi}. \tag*{\sqr53}\]

We conclude this section with the characterization of $\epsilon$-uniform mixing that follows directly from Theorem \ref{pgst_quotient} on regular non-bipartite graphs.

\begin{theorem}\label{thm:eum}
	Let $X$ be a $k$-regular non-bipartite graph on $n$ vertices. Let
	\[A = \sum_r k\cos\theta_r E_r\]
	be the spectral decomposition of the adjacency matrix of $X$, assuming $\theta_0=0$. Let $a$ be a vertex in $X$. Then $X$ admits local $\epsilon$-uniform mixing if and only if there is a vector $y\in\{\pm 1/\sqrt{kn}\}$ for which all of the following hold.
	\begin{enumerate}[(i)]
		\item $y$ lies in the intersection of $\col(I+R)+\col(D_t^T)$ and $\col(I-R) + \col(D_t^T)$.
		\item There exist real numbers $\delta_r$ where $\delta_0\in \{0, \pi\}$ such  that
		\begin{align*}
			D_t y &=\sqrt{k} \sum_{r} \cos\delta_r E_re_a,  \\
			D_h y &= \sqrt{k}  \sum_{r} \cos(\delta_r + \theta_r) E_re_a
		\end{align*}
		\item For any set $\{\ell_r: r\in \Phi_a\}$ of integers such that
		\begin{equation}\label{weak_eum}
			\sum_{r\in \Phi_a} \ell_r \theta_r =0,
		\end{equation}
		we have
		\[\sum_{r\in \Phi_a} \ell_r\delta_r \equiv 0 \pmod{2\pi}.\]
	\end{enumerate}
		Moreover, if (i) - (iii) hold with Equation \eqref{weak_eum} replaced by
	\[\sum_{r\in \Phi_a} \ell_r \theta_r \equiv 0 \pmod{2\pi},\]
	then $y$ can be approximated in integer steps.
\end{theorem}

\section{Equitable partitions \label{sec: equitable}}
A parition $\pi = \{\seq{C}{0}{1}{d}\}$ of the vertex set of a digraph $Y$ is \textsl{equitable} if the number of out-neighbors in $C_j$ of a vertex in $C_i$ depends only on $i$ and $j$. Alternatively, if $S$ is the characteristic matrix of $\pi$, then $\pi$ is equitable if there is a square matrix $W$ such that 
\[A(Y) S = SW;\]
this matrix $W$ is called the \textsl{quotient matrix}, and $W_{C_i, C_j}$ is precisely the number of out-neighbors in $C_i$ of any vertex in $C_i$. It is sometimes more convenient to work with the \textsl{normalized characteristic matrix} $\widehat{S}$, which is obtained from $S$ by scaling each of its columns to a unit vector:
\[\widehat{S} = (S^TS)^{-1/2} S.\]

The \textsl{line digraph} of a graph $X$, denoted $LD(X)$, is the digraph with the arcs of $X$ as vertices, where $(a,b)$ is adjacent to $(c,d)$ if $b=c$. Every equitable partition of $X$ naturally induces an equitable partition of $LD(X)$.

\begin{lemma}\label{lem: induced}
	 Let $\pi=\{\seq{C}{0}{1}{d}\}$ be an equitable partition of $X$. Let
	 \[C'_{ij} = \{(v_i, v_j): v_i\in C_i, v_j\in C_j, v_i \sim v_j\}.\]
	 Then  $\pi'=\{C'_{ij}\}_{i,j}$ is an equitable partition of $LD(X)$. 
\end{lemma}
\proof
Let $(v_i, v_j)$ be an arc of $X$ in the cell $C'_{ij}$. Let $C'_{\ell m}$ be any cell of $\pi'$. The number of neighbors of $(v_i, v_j)$ in $C'_{\ell m}$ is the size of the set
\[\{(v_{\ell}, v_m)\in C'_{\ell m}: v_{\ell} = v_j \}.\]
If $j\ne \ell$, then $C'_{\ell m}$ contains no arc with tail $v_j$, and so the above set is empty. If $j=\ell $, then the above set becomes
\[\{(v_j, v_m): v_m \in C_m  \text{ and }  v_m \sim v_j\},\]
whose size equals the number of neighbors in $C_m$ of $v_j$. In both cases, the number of neighbors of $(v_i, v_j)$ in $C'_{\ell m}$ depends only on $i,j,\ell, m$.
\qed

We will refer to the partition $\pi'$ defined in Lemma \ref{lem: induced} as the partition of the arcs of $X$ \textsl{induced by} $\pi$. The relation between $\pi$ and $\pi'$ can also be described algebraically.

\begin{lemma}
	Let $\pi$ be a partition of the vertices of $X$, with characteristic matrix $S$. Let $\pi'$ be the induced partition of the arcs of $X$, with characteristic matrix $S'$. Then $S'$ is a submatrix of $S\otimes S$.
\end{lemma}
\proof 
Let $(u,v)$ be an arc of $X$. Let $C'_{ij}$ be a cell of the induced partition. Then
\[S'_{(u,v), C'_{ij}} = \begin{cases}
	1,\quad \text{ if } u\in C_i, v\in C_j\\
	0,\quad \text{ otherwise,}
\end{cases}
\]
and 
\[(S\otimes S)_{(u,v), (C_i, C_j)} = S_{u,C_i} S_{v, C_j} =  \begin{cases}
	1,\quad \text{ if } u\in C_i, v\in C_j\\
	0,\quad \text{ otherwise.} \tag*{\sqr53}
\end{cases}\]

Recall that for a $k$-regular graph, the transition matrix of the quantum walk is
\[U = R\left(\frac{2}{k} D_t^T D_t -I\right).\]
The next result shows that the column space of the characteristc matrix of $\pi'$ is $U$-invariant, provided that $\pi'$ is induced by an equitable partition of the vertices of $X$.

\begin{lemma}
	Let $X$ be a regular graph. Let $\pi$ be an equitable partition of the vertices of $X$. Let $S$ be the characteristic matrix of the partition $\pi'$ induced by $\pi$. Then $\col(S)$ is $U$-invariant. 
\end{lemma}
\proof 
Notice that $D_h^TD_t=A(LD(X))$. So we may rewrite $U$ as
\[U =\frac{2}{k} D_h^T D_t -R =  \frac{2}{k} A(LD(X)) - R.\]
By Lemma \ref{lem: induced}, $\col(S)$ is $A(LD(X))$-invariant. Hence, it suffices to show that $\col(S)$ is $R$-invariant. Indeed, for each cell $C'_{ij}$ of $\pi'$, the image of the characteristic vector of $C'_{ij}$ under $R$ is the characteristic vector of $C'_{ji}$. Hence $\col(S)$ is $R$-invariant.
\qed

Finally, we prove an intuitive result about pretty good state transfer in any quantum walk that starts in a $U$-invariant subspace. Let $Q$ be the matrix whose columns form an orthonormal basis for this $U$-invariant subspace. Since $U$ is unitary, the orthogonal complement $(\col(Q))^{\perp}$ is also $U$-invariant. Hence $QQ^*$ commutes with $U$, and 
\[(Q^*UQ)^*(Q^*UQ) = Q^*U^*QQ^*UQ = Q^*QQ^*U^*UQ=Q^*U^*UQ \]
It follows that $Q^*UQ$ is unitary. In the case where $Q$ is the normalized characteristic matrix of a partition of the arcs induced by some equitable partition $\pi$ of the vertices, we refer to $Q^*UQ$ as the \textsl{quotient transition matrix} relative to $\pi$.

\begin{lemma}\label{pgst_quotient}
Let $Q$ be the matrix whose columns form an orthonormal basis for a $U$-invariant subspace. Then $U$ admits pretty good state transfer from $x=Qw$ to $y$ if and only if 
\begin{enumerate}[(i)]
	\item $y=Qz$ for some $z$, and 
	\item $Q^*UQ$ admits pretty good state transfer from $w$ to $z$.
\end{enumerate} 
\end{lemma}
\proof
Since $QQ^*$ commutes with $U$, it commutes with every eigenprojection of $U$. Thus for any time $t\in \re$,
\[U^t Qw = U^tQQ^*Qw = QQ^*U^tQw,\]
and so $U^tx$ lies in $\col(Q)$. As $y$ is a scalar multiple of some limit point of $\{U^t x\}_{t\in \re}$, it must also lie in $\col(Q)$. This shows that (i) is a necessary condition for pretty good state transfer to occur from $x=Qw$ to $y$.

Now let $y=Qz$ for some $z$. Using the spectral decomposition of $U$, we see that
\[(Q^*UQ)^t=Q^*U^t Q,\]
and so
\[Q(Q^*UQ)^t=QQ^*U^tQ=U^tQQ^*Q=U^tQ.\]
Thus
\[\norm{(Q^*UQ)^tw-z}=\norm{Q(Q^*UQ)^tw-Qz}= \norm{U^tQw-z}.\]
Therefore, $U$ admits pretty good state transfer from $x=Qw$ to $y=Qz$ if and only if (ii) holds. 
\qed

\section{Association schemes \label{sec: scheme}}

In this section, we study local $\epsilon$-uniform mixing on non-bipartite graphs in association schemes. 

An \textsl{association scheme} with $d$ classes is a set of $01$-matrices 
\[\mathcal{A}=\{\seq{A}{0}{1}{d}\}\]
that satisfies the following properties.
\begin{enumerate}[(i)]
	\item $A_0=I$.
	\item $A_0+A_1+\cdots+A_d=J$.
	\item $A_j^T\in\mathcal{A}$ for each $j$.
	\item $A_jA_{\ell} = A_{\ell}A_j\in\mathrm{span}(\mathcal{A})$.
\end{enumerate}
The algebra generated by $\mathcal{A}$ is called the \textsl{Bose-Mesner} algebra. It has two bases: $\{\seq{A}{0}{1}{d}\}$, whose elements are Schur idempotents, and $\{\seq{E}{0}{1}{d}\}$, whose elements are matrix idempotents. Thus there are scalars $p_r(s)$ such that
\[A_r = \sum_{s=0}^d p_r(s) E_s.\]
The change-of-basis matrix $P=(p_r(s))_{rs}$ is called the \textsl{matrix of eigenvalues} of $\mathcal{A}$. There are also non-negative integers $p_{i,j}(\ell)$ such that
\[A_i A_j = \sum_{\ell=0}^d p_{i,j}(\ell)A_{\ell}.\]
These $p_{i,j}(\ell)$'s are called the \textsl{intersection numbers} of $\mathcal{A}$. If all matrices in $\mathcal{A}$ are symmetric, then $\mathcal{A}$ is said to be \textsl{symmetric}, and each $A_j$ is the adjacency matrix of some graph. For more background on association schemes, see \cite[Ch 12]{Godsil1993}.

Given a symmetric association scheme $\mathcal{A}=\{\seq{A}{0}{1}{d}\}$, let $X$ be the graph with adjacency matrix $A_1$. For any vertex $a$ in $X$, the partition $\pi=\{\seq{C}{0}{1}{d}\}$ where 
\[C_j=\{u: (A_j)_{au}=1\}\]
is an equitable partition, with characteristic matrix
\[S=\pmat{A_0e_a & A_1 e_a & \cdots & A_d e_a}.\]
Let  $\pi'=\{C'_{ij}\}_{i,j}$ be the induced partition of the arcs, with characteristic matrix $S'$. As $C'_{01}$ is the set of outgoing arcs of $a$, the initial state $x_a$ is a scalar multiple of the corresponding column of $S'$, and so pretty good state transfer can only occur from $x_a$ to some state in $\col(S')$. Now, assuming $X$ has valency $k$, the transition matrix $U$ satisfies
\[US' = S'W,\]
where
\[W_{C'_{ij}, C'_{\ell m}} = \begin{cases}
	\frac{2}{k}-1,\quad \text{ if } i=m, j=\ell\\
	\frac{2}{k},\quad \text{ if } i\ne m, j=\ell\\
	0,\quad \text{ otherwise}.
\end{cases} \]
Clearly, $W$ does not depend on our choice of $a$, and neither does the quotient transition matrix. Applying Corollary \ref{pgst_quotient} with $Q=\widehat{S'}$ shows that local $\epsilon$-uniform mixing in symmetric association schemes implies $\epsilon$-uniform mixing.

\begin{theorem}
		Let $X$ be a graph in a symmetric association scheme. If $X$ admits local $\epsilon$-uniform mixing, then it admits $\epsilon$-uniform mixing.
\end{theorem}

Finally, using Lemma \ref{in-Dt} and Theorem \ref{thm:eum}, we find a connection between a special form of $\epsilon$-uniform mixing and real Hadamard matrices in Bose-Mesner algebras. We say a graph admits \textsl{local $\epsilon$-uniform mixing at $a$ towards $\col(D_t^T)$} if there is pretty good state transfer from $x_a$ to some state in $\col(D_t^T)$, and if this happens for every vertex $a$, we say $X$ admits \textsl{$\epsilon$-uniform mixing towards $\col(D_t^t)$}. We will show in the next section that, up to multiplication by $R$, this is the only form of $\epsilon$-uniform mixing that could occur on a complete graph or a strongly regular graph.

\begin{theorem}\label{eum_had}
	Let $X$ be a non-bipartite graph in a symmetric association scheme, and let
	\[A = \sum_r k\cos\theta_r E_r\]
	be the spectral decomposition of the adjacency matrix of $X$, assuming $\theta_0=0$. Then $X$ admits $\epsilon$-uniform mixing towards $\col(D_t^T)$ if and only if the following hold.
	\begin{enumerate}[(i)]
		\item The Bose-Mesner algebra contains a real Hadamard matrix $H$.
		\item Assuming, by possibly negating $H$, that
		\[H =\sqrt{n}\left( E_0 + \sum_{r\ne 0} (-1)^{\sigma_r} E_r\right).\]
		 For any set $\{\ell_r: r\in \Phi_a-\{r\}\}$ of integers such that
		\begin{equation}\label{weak_eum_Dt}
			\sum_{r\in \Phi_a-\{0\}} \ell_r \theta_r =0,
		\end{equation}
		the integer
		\[\sum_{r\in \Phi_a-\{0\}} \ell_r\sigma_r \]
		is even.
	\end{enumerate}
	Moreover, if (i) - (iii) hold with Equation \eqref{weak_eum_Dt} replaced by
	\[\sum_{r\in \Phi_a-\{0\}} \ell_r \theta_r \equiv 0 \pmod{2\pi},\]
	then the target state can be approximated in integer steps.
\end{theorem}
\proof
Suppose pretty good state transfer occurs from $x_a$ to $D_t^Tw$ for some real flat $w$. Then $x_a$ is strongly cospectral to $D_t^Tw$, and By Lemma \ref{in-Dt}, there are integers $\sigma_r$ such that
\[
w =\frac{1}{\sqrt{k}}\sum_{r} (-1)^{\sigma_r} E_re_a.
\]
Since $A$ lies in an association scheme $\{\seq{A}{0}{1}{d}\}$, each $E_r$ is a linear combination of $\seq{A}{0}{1}{d}$, and so 
\[\sqrt{n}\sum_{r} (-1)^{\sigma_r} E_r\]
has flat entries and $\pm\sqrt{n}$ eigenvalues; that is, it is a real Hadamard matrix. The rest of the result follows from Theorem \ref{thm:eum}.
\qed

\section{Distance regular graphs}
A graph is \textsl{distance regular} if for any two vertices $u$ and $v$ at distance $\ell$, the number of vertices at distance $i$ from $u$ and at distance $j$ from $v$ depends only on $i, k, \ell$. Given a distance regular graph with diameter $d$, its distance matrices $\seq{A}{0}{1}{d}$ form an association scheme with $d$ classes, and the intersection numbers are denoted $a_i$, $b_i$ and $c_i$, where 
\[a_i = p_{1, i}(i),\quad b_i = p_{1,i+1}(i),\quad c_i = p_{1, i-1}(i).\]
By our discussion in Sections \ref{sec: equitable} and \ref{sec: scheme}, if a distance regular graph $X$ admits pretty good state transfer from $x_a$ to $y$, then $y$ is constant on the cells 
\[C_{ij} = \{(u, v): u\sim v, \dist{a}{u}=i, \dist{a}{v}=j\}\]
where $j=i-1$, $i$ or $i+1$. Let $S$ be the characteristic matrix of this parititon. We will investigate real states of the form $y=Sw$ that are strongly cospectral to $x_a$.

Our first observation is that any such state lies in $\col(I+R) + \col(D_t^T)$.

\begin{lemma}
	If $X$ is a distance regular graph, then $\col(S)$ is a subspace of $\col(I+R) + \col(D_t^T)$.
\end{lemma}
\proof
We show by induction that the characteristic vectors of $C_{i,j}$ lies in $\col(I+R)+\col(D_t^T)$. 

First note that the cell $C_{0,1}$  consists of all outgoing arcs of $a$, so its characteristic vector $D_t^Te_a$ lies in $\col(D_t^T)$. As $\col(I+R)+\col(D_t^T)$ is $R$-invariant, it also contains the characteristic vector of $C_{1,0}$. The third cell $C_{1,1}$ consists of all arcs between the neighbors of $a$, if any, and its characteristic vector lies in $\col(I+R)$.

Now suppose for any $i,j, \le \ell$, the characteristic vector of $C_{i,j}$ stays in $\col(I+R)+\col(D_t^T)$. Consider $C_{\ell, \ell+1}$. We may write its characteristic vector as 
\[S e_{\ell, \ell+1} = \sum_{u: \dist{a}{u}=\ell} D_t^T e_u - S e_{\ell, \ell-1} - S e_{\ell, \ell},\]
which, by the induction hypothesis, lies in $\col(I+R)+\col(D_t^T)$. As $C_{\ell+1, \ell}$ consists of the reversed arcs, its characteristic also lies in $\col(I+R)+\col(D_t^T)$. Finally, the characteristic vector of $C_{\ell, \ell}$ is clearly in $\col(I+R)$.
\qed

Recall that $B$ denotes the unsigned vertex-edge incidence matrix, and $M$ denotes the unsigned arc-edge incidence matrix. Let $K$ be the matrix whose rows are indexed by the edges of $X$ and columns by the cells $C_{ij}$, where
\[K_{ \{u,v\}, C_{i,j}}=\begin{cases}
	2,\quad \text{if $i=j$ and $\dist{a}{u}=\dist{a}{v}=i$}\\
	1,\quad \text{if $i\ne j$,  $\dist{a}{u}=i$ and $\dist{a}{v}=j$}\\
	0,\quad \text{otherwise}.
\end{cases}\]
It is straightforward to verify that $K=M^TS$. This gives an alternative, and in some cases more convenient, description of states of the form $Sw$ that lie in $\col(I-R)+\col(D_t^T)$.

\begin{lemma}
		If $X$ is a distance regular graph, then $Sw$ lies in $\col(I-R)+\col(D_t^T)$  if and only if $Kw$ lies in $\col(B^T)$.
\end{lemma}
\proof
 By Lemma \ref{props}, the orthogonal complement of $\col(I-R)+\col(D_t^T)$ is $M\ker(B)$. Hence, $Sw$ lies in $\col(I-R)+\col(D_t^T)$ if and only if it is orthogonal to $M\ker(B)$, that is, $Kw$ lies in $\col(B^T)$.
 \qed

With the above observations, we now reduce Theorem \ref{strcosp} to a characterization involving the parameters of the association scheme. For ease of notation, let $R'$ be the permutation that swaps cell $C_{i,j}$ with cell $C_{j,i}$.

\begin{theorem}\label{eum_drg}
	Let $X$ be a non-bipartite distance regular graph on $n$ vertices of diameter $d$. Let $\mathcal{A}$ be the corresponding association scheme, with matrix of eigenvalues $P=(p_r(s))_{rs}$. For $r=0,1,\cdots,d$, let $\theta_r = p_1(r)$. Then $x_a$ is strongly cospectral to a real state $y$ if and only if $y=Sw$ for some $w$ satisfying the following.
	\begin{enumerate}[(i)]
		\item $Kw$ lies in $\col(B^T)$.
		\item There are scalars $\delta_r$ where $\delta_0\in \{0, \pi\}$ such that 
			\[\sqrt{k} \pmat{ \cos\delta_0 \\\cos\delta_1 \\ \vdots \\ \cos\delta_d}=P\pmat{b_0 &  &  & &&&&&  \\ &c_1 & a_1 & b_1&&&& \\ & & && &\ddots &&&\\  & &&  & &&&c_d & a_d}w\]
			and 
				\[\sqrt{k} \pmat{\cos(\delta_0+\theta_0) \\ \cos(\delta_1+\theta_1) \\ \vdots \\ \cos(\delta_d+\theta_d)}=P\pmat{b_0 &  &  & &&&&&  \\ &c_1 & a_1 & b_1&&&& \\ & & && &\ddots &&&\\  & &&  & &&&c_d & a_d}R'w\]
	\end{enumerate}
\end{theorem}
\proof
It suffices to show that Condition (ii) is equivalent to Theorem \ref{strcosp} (ii). Let $\{\seq{E}{0}{1}{d}\}$ be the idempotents of $\mathcal{A}$. With $y=Sw$, the first identity in Theorem \ref{strcosp} (ii) becomes  
\begin{align}
	D_t Sw &=\sqrt{k}\sum_{r} \cos\delta_r E_r e_a \label{eqn: E_Dt}\\
	&=\sqrt{k} \pmat{E_0 e_a & \cdots & E_d e_a} \pmat{\cos\delta_0 \\ \cos\delta_1 \\ \vdots \\ \cos\delta_d}\notag\\
	&=\sqrt{k} \pmat{A_0 e_a & \cdots & A_d e_a} P^{-1}\pmat{\cos\delta_0 \\ \cos\delta_1 \\ \vdots \\ \cos\delta_d}. \notag 
\end{align}
Note that $\Delta= \pmat{A_0 e_a & \cdots & A_d e_a}$ is the characteristic matrix of the distance partition of $V(X)$ relative to $a$, and $\Delta^T\Delta$ is a diagonal matrix whose $ii$-entry counts the vertices at distance $i$ from $a$. Moreover, 
\[(\Delta^T D_t)_{i, (u, v)}= 
\begin{cases}
	1,\quad \text{if $\dist{a}{u}=i$}\\
	0,\quad \text{otherwise}.
\end{cases}\]
Hence, $(\Delta^T D_t S)_{\ell, C_{i,j}}$ counts the arcs in cell $C_{i,j}$ whose tails are at distance $\ell$ from $a$. We then have
\[(\Delta^T\Delta)^{-1}\Delta^T D_t S =\pmat{b_0 &  &  & &&&&&  \\ &c_1 & a_1 & b_1&&&& \\ & & && &\ddots &&&\\  & &&  & &&&c_d & a_d}\]
By Equation \eqref{eqn: E_Dt}, 
\begin{equation}\label{eqn:simple_Dt}
	\sqrt{k} \pmat{\cos\delta_0 \\\cos\delta_1 \\ \vdots \\ \cos\delta_d}=P\pmat{b_0 &  &  & &&&&&  \\ &c_1 & a_1 & b_1&&&& \\ & & && &\ddots &&&\\  & &&  & &&&c_d & a_d}w.
	\end{equation}
Conversely, if Equation \eqref{eqn:simple_Dt} holds, then we may multiply both sides by $\pmat{E_0 e_a & \cdots & E_d e_a}$ to derive Equation \eqref{eqn: E_Dt}:
\begin{align*}
	\sqrt{k} \pmat{E_0 e_a & \cdots & E_d e_a} \pmat{\cos\delta_0 \\ \cos\delta_1 \\ \vdots \\ \cos\delta_d} &= \Delta\pmat{b_0 &  &  & &&&&&  \\ &c_1 & a_1 & b_1&&&& \\ & & && &\ddots &&&\\  & &&  & &&&c_d & a_d} w\\
	&=D_t S w.
\end{align*}
A similar argument shows that the second identity in Theorem \ref{strcosp} (ii) is equivalent to 
\[\sqrt{k} \pmat{\cos(\delta_0+\theta_0) \\ \cos(\delta_1+\theta_1) \\ \vdots \\ \cos(\delta_d+\theta_d)}=P\pmat{b_0 &  &  & &&&&&  \\ &c_1 & a_1 & b_1&&&& \\ & & && &\ddots &&&\\  & &&  & &&&c_d & a_d}R'w\tag*{\sqr53}\]

This immediately gives a classification of complete graphs that admits local $\epsilon$-uniform mixing.

\begin{theorem}
	The only complete graphs that admit local $\epsilon$-uniform mixing are $K_2$ and $K_4$.
\end{theorem}
\proof
Suppose $K_n$, where $n\ge 3$, admits local $\epsilon$-uniform mixing at vertex $a$. Then there is pretty good state transfer from $x_a$ to $y=\frac{1}{\sqrt{n(n-1)}}Sz$ for some $z\in \{\pm1\}^3$. By Lemma \ref{sc_R}, $x_a$ is strongly cospectral to both $y$ and $Ry$.

If $z_{10} = z_{11}$ or $(R'z)_{10}=(R'z)_{01}$, then one of $y$ and $Ry$ lies in $\col(D_t^T)$, and it follows from Theorem \ref{eum_had} that the adjacency algebra of $K_n$ contains a real Hadamard matrix. This restricts $n$ to $4$.

Otherwise, we may assume $z_{01} = z_{10} = - z_{11}=1$. We have
\[b_0 = n-1,\quad c_1=1, \quad a_1 = n-2\]
and
\[P = \pmat{1 & n-1\\1 & -1}.\]
So by Theorem \ref{eum_drg} (ii), there exists $\delta_1\in \re$ such that $\cos\delta_1=\cos(\delta_1+\theta_1)$ and 
\[(n-1)\sqrt{n} \pmat{\pm1 \\ \cos\delta_1}= \pmat{1 & n-1\\1 & -1}\pmat{n-1 & & \\ & 1 & -1} \pmat{1 \\ 1\\ -1}\]
However, with $\cos\theta_1 = -1/(n-1)$, this implies $n=9$ and 
\[\cos\delta_1 = \frac{1}{4},\quad \cos\theta_1 = -\frac{1}{8},\]
which do not satisfy $\cos\delta_1 = \cos(\delta_1+\theta_1)$. Therefore, $y$ is not strongly cospectral to $x_a$, a contradiction.
\qed

\section{Strongly regular graphs}

A \textsl{strongly regular graph} is a distance regular graph of diameter $2$. Its parameters are determined by the number of vertices $n$, the valency $k$, the number $a_1$ of common neighbors of any pair of adjacent vertices, and the number $c_2$ of common neighbors of any pair of non-adjacent vertices. These parameters satisfy the relation
\[k(k-a_1-1) = (n-k-1)c_2.\]
Let $k>\lambda_1>\lambda_2$ be the eigenvalues of an $(n,k,a_1,c_2)$-strongly regular graph. 
Then
\[\lambda_r = \frac{a_1-c_2 \pm \sqrt{ (a_1-c_2)^2 + 4(k-c_2)}}{2},\]

In this section, we determine all $(n,k,a_1, c_2)$-strongly regular graphs that admit $\epsilon$-uniform mixing. It is shown in Lemma \ref{bip_srg} that this does not occur on $K_{k,k}$ for $k\ge 2$. Hence, we will consider non-bipartite strongly regular graphs only. By Theorem \ref{eum_had}, if any of these graphs admits $\epsilon$-uniform mixing towards $\col(D_t^T)$, then its adjacency algebra contains a real Hadamard matrix. Such strongly regular graphs have been fully characterized and constructed.

\begin{theorem}\cite{Goethals1970, Haemers2010,Brouwer1984}\label{reghad}
	The adjacency algebra of a strongly regular graph $X$ contains a real Hadamard matrix if and only if $X$ or $\comp{X}$ has parameters $(4m^2, 2m^2\pm m, m^2\pm m, m^2\pm m)$.
\end{theorem}

Using Theorem \ref{eum_had}, we show all these graphs except for $C_4$ and $2K_2$ admit $\epsilon$-uniform mixing. We need the following two results about \textsl{pure geodetic angles}, that is, angles whose six squared trigonometric functions are either rational or infinite.

\begin{theorem}\cite[Ch 11]{Bergen2009}\label{thm:tan_rat_pi} %
	If $q\in \rats$ and $q\pi$ is a pure geodetic angle, then 
	\[\tan(q\pi) \in \left\{0, \pm{\sqrt{3}}, \pm \frac{1}{\sqrt{3}}, \pm 1\right\}.\]
\end{theorem}

\begin{theorem}\cite{Conway1999}\label{thm:rat_lb_geo}
	If the value of a rational linear combination of pure geodetic angles is a rational multiple of $\pi$, then so is the value of its restriction to those angles whose tangents are rational multiple of any given square root.
\end{theorem}

We get four infinite families of graphs that admit $\epsilon$-uniform mixing towards $\col(D_t^T)$, in integer steps.

\begin{theorem}\label{egs}
	If $X$ or $\comp{X}$ is a $(4m^2, 2m^2\pm m, m^2\pm m, m^2\pm m)$-strongly regular graph, where $m\ge 2$, then $X$ admits $\epsilon$-uniform mixing in integer steps.
\end{theorem}
\proof
We prove this for the parameter set $(4m^2, 2m^2+ m, m^2+ m, m^2+ m)$; the other case follows similarly. Let $X$ be a $(4m^2, 2m^2+ m, m^2+ m, m^2+ m)$-strongly regular graph. By Theorem \ref{reghad}, its adjacency algebra contains a real Hadamard matrix. Note that the eigenvalues of $X$ are 
\[k=2m^2+m,\quad \lambda_1=m,\quad  \lambda_2=-m,\]
and all of them lie in the eigenvalue support of every vertex. 
Let 
\[\theta_1 =\arccos\left(\frac{\lambda_1}{k}\right),\quad \theta_2 =\arccos\left(\frac{\lambda_2}{k}\right),\]
and let $\ell_1$ and $\ell_2$ be integers such that 
\[\ell_1\theta_1+\ell_2\theta_2\equiv 0\pmod{2\pi}.\]
Since $\lambda_2=-\lambda_1$, we have $\theta_2 = \pi-\theta_1$ and 
\[\ell_2 \pi + (\ell_1-\ell_2)\theta_1\equiv 0 \pmod{2\pi}.\]
Moreover, by Theorem \ref{thm:tan_rat_pi}, the angle $\theta_1 = \arccos\left(\frac{1}{2m+1}\right)$ is not a rational multiple of $\pi$ for any positive integer $m$. Hence $\ell_1-\ell_2=0$, and both $\ell_1$ and $\ell_2$ must be even. By Theorem \ref{eum_had}, $X$ admits $\epsilon$-uniform mixing in integer steps. 

Now consider the complement $\comp{X}$. It has parameters $(4m^2, 2m^2-m-1, m^2-m-2, m^2-m)$, and its eigenvalues are 
\[k = 2m^2-m-1,\quad \lambda_1 = m-1,\quad \lambda_2=-m-1.\]
Again, all of them lie in the eigenvalue support of every vertex. 
Let 
\[\theta_1 =\arccos\left(\frac{\lambda_1}{k}\right),\quad \theta_2 =\arccos\left(\frac{\lambda_2}{k}\right),\]
and let $\ell_1$ and $\ell_2$ be integers such that 
\[\ell_1\theta_1+\ell_2\theta_2\equiv 0 \pmod{2\pi}.\]
We have
\[\tan\theta_1 = 2\sqrt{m(m+1)},\quad \tan\theta_2 = -\frac{2m}{m+1} \sqrt{m^2-m-1}.\]
If $\tan\theta_1$ and $\tan\theta_2$ were both rational multiples of $\sqrt{\Delta}$ for some square-free integer $\Delta$, then 
\[m(m+1)(m^2-m-1)=m^4 - 2m^2 - m.\] 
would be a perfect square. However, this is impossible for $m\ge 2$ as it lies strictly between $(m^2-2)^2$ and $(m^2-1)^2$. Hence, by Theorem \ref{thm:rat_lb_geo}, both $\ell_1\theta_1$ and $\ell_2\theta_2$ are rational multiples of $\pi$. As $\tan\theta_1, \tan\theta_2 \notin \{0,\pm\sqrt{3}, \pm\frac{1}{\sqrt{3}}, \pm1\}$, both $\ell_1$ and $\ell_2$ are zero. Therefore by Theorem \ref{eum_had}, $\comp{X}$ admits $\epsilon$-uniform mixing in integer steps.
\qed

Our next goal is to show that this is the only form of $\epsilon$-uniform mixing that could occur on a strongly regular graph. Given the eigenvalue $k>\lambda_1>\lambda_2$ of a strongly regular graph, the matrix of eigenvalues of its association scheme is
\[P = \pmat{
	1 & k & n-k-1\\
	1 & \lambda_1 & -1-\lambda_1\\
	1 & \lambda_2 & -1-\lambda_2}.\]

\begin{lemma}\label{feas}
	Let $X$ be an $(n,k,a_1,c_2)$-strongly regular graph with $c_2\ne k$. Suppose $X$ admits $\epsilon$-uniform mixing. Then $n$ is a perfect square, and there exists a vector
	\[z = \pmat{z_{01} & z_{10} & z_{11} & z_{12} & z_{21} & z_{22}}^T\in \{\pm1\}^6\]
	(with the understanding that $z_{11}$ is not present if $a_1=0)$ that satisfies the following.
	\begin{enumerate}[(i)]
		\item $Kz$ lies in $\col(B^T)$.
		\item $\pm \sqrt{n} = z_{01} + z_{10} + a_1 z_{11} + (k-a_1-1)(z_{12}+z_{21}-z_{22}) + (n-k-1) z_{22}$.
		\item For each  eigenvalue $\lambda_r\ne k$ of $X$,
	\begin{align*}
		&2kn(k^2-\lambda_r^2)\\
		=& (k-\lambda_r)((k+\lambda_r)(z_{01}+z_{10})+2\lambda_r a_1 z_{11} + (\lambda_r b_1 - (\lambda_r+1)c_2)(z_{12}+z_{21})-2(\lambda_r+1)a_2 z_{22})^2\\
		&+(k+\lambda_r)((k-\lambda_r)(z_{01}-z_{10})+ (\lambda_r b_1 + (\lambda_r+1)c_2)(z_{12}-z_{21}))^2.
	\end{align*}
		\end{enumerate}
\end{lemma}
\proof
Corollary \ref{perfect_square} implies that $n$ is a perfect square. Theorem \ref{eum_drg} (i) implies that $Kz$ lies in $\col(B^T)$. For $r=1,2$, write $\lambda_r = k\cos\theta_r$. By Theorem \ref{eum_drg} (ii), there exist real numbers $\delta_1, \delta_2$ such that
\begin{align}
	k\sqrt{n}\pmat{\pm1 \\ \cos\delta_1 \\ \cos\delta_2} &= P\pmat{k&&&&&\\&1&a_1&b_1&&&\\&&&&c_2 & a_2} z \label{tail}\\
	k\sqrt{n}\pmat{\pm1 \\ \cos(\delta_1+\theta_1) \\ \cos(\delta_2+\theta_2)} &= P\pmat{k&&&&&\\&1&a_1&b_1&&&\\&&&&c_2 & a_2} R'z \label{head}
\end{align} 
Comparing the first entry of both sides of Equation \eqref{tail} shows 
\[\pm k\sqrt{n} = k(z_{01}+z_{10}) + ka_1z_{11}+kb_1(z_{12}+z_{21}) +(n-k-1)a_2z_{22},\]
which, as $(n-k-1)c_2=kb_1$, is equivalent to the expression in (ii). Now compare the remaining entries of both sides of Equation \eqref{tail} and \eqref{head}. We see that for $r=1,2$,
\begin{align*}
	k\sqrt{n} \cos\delta_r &= kz_{01}+\lambda_r z_{10} + a_1\lambda_r z_{11} + b_1\lambda_r z_{12} - (1+\lambda_r) (c_2 z_{21}-(1+\lambda) a_2 z_{22}), \\
	k\sqrt{n} \cos(\delta_r+\theta_r) &= kz_{10}+\lambda_r z_{01} + a_1\lambda_r z_11 + b_1 \lambda_r z_{21} - (1+\lambda_r) c_2 z_{12} - (1+\lambda)a_2z_{22}.
\end{align*}
Taking the sum and difference of the above gives
\begin{align}
	&2k\sqrt{n} \cos\left(\delta_r +\frac{\theta_r}{2}\right)\cos\left(\frac{\theta_r}{2}\right)  \label{cos}\\
	&=(k+\lambda_r)(z_{10}+z_{01}) + 2a_1\lambda_r z_{11} + (b_1\lambda_r - c_2(1+\lambda_r))(z_{12}+z_{21}) - 2(\lambda_r+1) a_2 z_{22} \notag
\end{align}
and 
\begin{align}
	&2k\sqrt{n} \sin\left(\delta_r +\frac{\theta_r}{2}\right)\sin\left(\frac{\theta_r}{2}\right) \label{sin}\\ 
	&=(k-\lambda_r)(z_{10}-z_{01}) + (b_1\lambda_r + c_2(1+\lambda_r)) (z_{12}-z_{21}). \notag
\end{align}
Since $\cos \theta_r\ne \pm 1$ and 
\begin{equation}\label{doubleangle}
	\cos \theta_r = 2\cos^2 (\theta_r/2)-1 = 1-2\sin^2(\theta_r/2),
\end{equation}
neither $\cos (\theta_r/2)$ nor $\sin (\theta_r/2)$ is zero. Thus, we may rearrange Equations \eqref{cos} and Equation \eqref{sin} to get
\[4k^2n = \left(\frac{\eqref{cos}}{\cos(\theta_r/2)}\right)^2 + \left(\frac{\eqref{sin}}{\sin(\theta_r/2)}\right)^2.\]
Using Equation \ref{doubleangle} again and the fact that $\lambda_r = k\cos\theta_r$, we obtain the desired expression in (iii).
\qed

We first rule out $\epsilon$-uniform mixing on triangle-free strongly regular graphs that are not the Clebsch graph. The following bound due to Biggs \cite{Biggs2009} turns out to be helpful.

\begin{lemma}\cite{Biggs2009}
	An $(n,k,0,c_2)$-strongly regular graph that is not $C_5$, the Clebsch graph, the Hoffman-Singleton graph, the Sims-Gewirtz graph, the M22 graph, or the Higman-Sims graph, satisfies
	\[k\ge \frac{7}{2}c_2 + \frac{25}{4}.\]
\end{lemma}

\begin{theorem}
	The only triangle-free strongly regular graph that admits $\epsilon$-uniform mixing is the Clebsch graph.
\end{theorem}
\proof
Let $X$ be an $(n,k,0,c_2)$-strongly regular graph that admits $\epsilon$-uniform mixing. Then $X$ is not bipartite and $n$ is a perfect square. If $\epsilon$-uniform mixing occurs towards $\col(D_t^T)$, then by Theorem \ref{egs}, $X$ has parameters $(16,5,0,2)$, and the unique strongly regular graph with these parameters is the Clebsch graph. So assume $\epsilon$-uniform mixing occurs but not towards $\col(D_t^T)$. Then either $X$ is the Higman-Sims graph, or $k\ge 7c_2/2+25/4$.

By Lemma \ref{feas} (ii), there is a $\{\pm1\}$-vector $z$ such that
\begin{equation}\label{feas1_a=0}
	\pm \sqrt{n} = z_{01} + z_{10} + (k-1) (z_{12}+z_{21}) + (n-2k) z_{22}.
\end{equation}
Suppose first that $z_{12}=z_{21}$. If $z_{21}=z_{22}$, then to ensure $Sz, RSz\notin \col(D_t^T)$, we need $z_{01}=z_{10}=-z_{12}$, but this is impossible as it would reduce Equation \eqref{feas1_a=0} to $\pm \sqrt{n} = n-4$. Hence $z_{21}=-z_{22}$, and up to a sign and a swap between $z$ and $R'z$, $z$ is one of 
\[\pmat{1&1&1&1&-1}^T,\quad \pmat{1&-1&1&1&-1}^T,\quad \pmat{1&1&-1&-1&1}^T.\]
Plugging these into Equation \eqref{feas1_a=0}, we see that $k$ takes one of the following values:
\begin{enumerate}[(i)]
	\item $4k = n\pm \sqrt{n}$ if $z=\pmat{1&1&1&1&-1}^T$
	\item $4k=n\pm\sqrt{n}+2$ if $z= \pmat{1&-1&1&1&-1}^T$
	\item $4k=n\pm \sqrt{n}+4$ if $z=\pmat{1&1&-1&-1&1}^T$
\end{enumerate}
The parameters of the Higman-Sims graph, $(100, 22, 0, 6)$, satisfy $4k=n-\sqrt{n}-2$. However, with $z= \pmat{1&-1&1&1&-1}^T$ and eigenvalue $\lambda_1 = 2$, the equation in Lemma \ref{feas} (iii) does not hold. Hence $X$ must be a different strongly regular graph with $k\ge 7c_2/2+25/4$. Since $k(k-1)=(n-k-1)c_2$, the choices of $k$ in (i), (ii) and (iii) limit $n$ to perfect squares no greater than $324$, but unless $n\in \{ 16, 100\}$, no such graph exists according to Andries Brouwer's database.

Now suppose $z_{12}=-z_{21}$. By possibly replacing $z$ with $R'z$, we may assume that $z_{21}=z_{22}$. To ensure $Sz, RSz\notin \col(D_t^T)$, $z$ or $-z$  is one of
\[\pmat{1&-1&1&-1&-1}^T, \quad \pmat{1&1&-1&1&1}.\]
By Equation \eqref{feas1_a=0} , we obtain two possibilities for $k$:
\begin{enumerate}[(i)]
\item$2k=n\pm \sqrt{n}$ if $z=\pmat{1&-1&1&-1&-1}^T$
\item $2k=n\pm \sqrt{n}+2$ if $z=\pmat{1&1&-1&1&1}$
\end{enumerate}
Each of these, together with $k\ge 7c_2/2+25/4$ and $k(k-1)=(n-k-1)c_2$, restricts $n$ to perfect squares no greater than $4$. Hence the Clebsch graph is the only triangle-free strongly regular graph that admits $\epsilon$-uniform mixing.
\qed

We now move on to non-bipartite strongly regular graphs that contain a triangle. Surprisingly, the fact that $a_1\ne 0$ greatly narrows down the choices of $z$ that could satisfy Lemma \ref{feas}. 

Given a strongly regular graph and a vertex $a$, the first and second \textsl{subconstituents} relative to $a$, denoted $\Gamma_1(a)$ and $\Gamma_2(a)$, are the subgraphs of $X$ induced by the neighbors of $a$ and the non-neighbors of $a$, respectively. We prove a technical lemma utilizing certain even cycles across these subconstituents.
	
\begin{lemma}\label{cycle}
	Let $X$ be an $(n,k,a_1,c_2)$-strongly regular graphs with $a_1\ne 0$ and $c_2\ne k$. Let $a$ be any vertex. If $x_a$ is strongly cospectral to $Sw$ for some $w\in \re^6$, then
	\[w_{11}+w_{22} = w_{12}+w_{21}.\] 
\end{lemma}	
\proof
Let $u_2$ be a vertex in $\Gamma_2(a)$. First note that $u_2$ lies in a path $(u_1, u_2, v_2, v_1)$ where $u_1, v_1 \in \Gamma_1(a)$, $v_2\in \Gamma_2(a)$ and $v_2\sim v_1$. If not, then $c_2=1$ and some vertex in $\Gamma_1(a)$ would be adjacent to $u_2$ and all its neighbors in $\Gamma_2(a)$, which would imply
\[b_1\ge a_2 + 1=a_2+c_2 = k,\]
a contradiction. Now if $u_1\sim v_1$, then $(u_1, u_2, v_2, v_1, u_1)$ forms a $4$-cycle, which we denote by $C_1$. Otherwise, $v_2$ is adjacent to some vertex $w$ in $\Gamma_1(a)$, and $(u_1, u_2, v_2, v_1, w, a, u_1)$ forms a $6$-cycle, which we denote by $C_2$.

By Lemma \ref{feas} (i), $Kw$ lies in $\col(B^T)$, so it is orthogonal to $\ker(B)$. For every even cycle in $X$, the signed characteristic vector of its edges, where the sign alternates along the cycle, lies in $\ker(B)$. Thus, the existence of $C_1$ implies
\[(w_{12}+w_{21}) - 2w_{22} + (w_{12}+w_{21}) - 2w_{11}=0,\]
and the existence of $C_2$ implies
\[(w_{12}+w_{21}) - 2w_{22} + (w_{12}+w_{21}) - 2w_{11} + (w_{01}+w_{10}) - (w_{01}+w_{10})=0.\]
Both reduce to $w_{11}+w_{22} = w_{12}+w_{21}$.
\qed

As a result of this lemma, it remains to check the feasibility conditions in Lemma \ref{feas} for two cases, and prove neither case works.

\begin{theorem}\label{simplefeas}
	Let $X$ be an $(n,k,a_1,c_2)$-strongly regular graphs with $a_1\ne 0$ and $c_2\ne k$. If $\epsilon$-uniform mixing occurs on $X$, then it occurs towards $\col(D_t^T)$.
\end{theorem}
\proof
Suppose pretty good state transfer occurs from $x_a$ to $y=\frac{1}{kn}Sz$ where $z\in\{\pm 1\}^6$, and for the sake of contradiction, neither $y$ nor $Ry$ lies in $\col(D_t^T)$. By Lemma \ref{cycle}, $z_{11}+z_{22} = z_{12}+z_{21}$, and so Lemma \ref{feas} (ii) implies
	\begin{equation}\label{feas1_a>0}
		\pm\sqrt{n} = z_{01} + z_{10} + (k-1)z_{11} + (n-k-1) z_{22}.
		\end{equation}
	 Suppose first that $z_{11}=z_{22}$. Then $z_{12}=z_{21}=z_{11}$. To ensure $y, Ry \notin \col(D_t^T)$, we need $z_{01}=z_{10}=-z_{11}=-z_{22}$, but this is impossible as it would reduce Equation \eqref{feas1_a>0} to $\pm\sqrt{n} = n-4$. Hence $z_{11}=-z_{22}$. It follows that $z_{12}=-z_{21}$. By replacing $z$ with $R'z$ if necessary, we may assume that $z_{11}=z_{12}$. Now to ensure $y, Ry \notin \col(D_t^T)$, $z$ or $-z$ must satisfy one of 
	 \[\pmat{1&-1&1&1&-1&-1}^T,\quad \pmat{1&1&-1&-1&1&1}^T.\]
	 Plugging the values of $z_{01}, z_{10}, z_{11},z_{22}$ into Equation \eqref{feas1_a>0} shows
	 \begin{enumerate}[(i)]
	 	\item$2k=n\pm \sqrt{n}$ if $z=\pmat{1&-1&1&1&-1&-1}^T$, and
	 	\item $2k=n\pm \sqrt{n}+2$ if $z=\pmat{1&1&-1&-1&1&1}^T$.
	 \end{enumerate}
To see if these choices meet the constraint in Lemma \ref{feas} (iii), let
	\begin{align*}
	f_r
	=& (k-\lambda_r)((k+\lambda_r)(z_{01}+z_{10})+2\lambda_r a_1 z_{11} + (\lambda_r b_1 - (\lambda_r+1)c_2)(z_{12}+z_{21})-2(\lambda_r+1)a_2 z_{22})^2\\
	&+(k+\lambda_r)((k-\lambda_r)(z_{01}-z_{10})+ (\lambda_r b_1 + (\lambda_r+1)c_2)(z_{12}-z_{21}))^2-2kn(k^2-\lambda_r^2),
\end{align*}
and define $g=f_1+f_2$ and $h=f_1-f_2$. Note that $f_1=f_2=0$ if and only if $g=h=0$, and we can express
\begin{align*}
	\lambda_1+\lambda_2&=a_1-c_2,\\
	\lambda_1-\lambda_2&=\sqrt{(a_1-c_2)^2+4(k-c_2)}.
\end{align*}
Thus, with the above two choices of $k$ and $z$, the parameters of $X$ satisfy
\begin{enumerate}[(1)]
	\item $k(k-a_1-1)=(n-k-1)c_2$,
	\item $g=0$, and
	\item $h=0$.
\end{enumerate}
However, a computation using SageMath shows that these cannot simultaneously hold. To illustrate, let $\ell=\sqrt{n}$. For Case (i), these feasibility conditions imply 
\[\ell^4\pm \ell^3-(4b_1+4)\ell^2\mp(4b_1+4)\ell + 16b_1=0,\]
from which it follows that $\ell$ divides $16b_1$ and so $\ell^2\pm\ell-4$ divides $16$. But then $\ell\in\{2,3,4,5\}$, none of which gives rise to a positive integer $b_1$. Similarly, for Case (i) in Corollary \ref{simplefeas}, these feasibility conditions imply
\[\ell^3\mp \ell^2-(4b_1+4)\ell \pm 4b_1=0,\]
from which it follows that $\ell$ divides $4b_1$ and so $\ell\mp 1$ divides $4$. The only feasible integer pair of $(\ell, b_1)$ is $(5,5)$, but it does not correspond to a strongly regular graph.
	\qed

This completes the classfication of strongly regular graphs that admit $\epsilon$-uniform mixing.

\begin{theorem}
	A strongly regular graph $X$ admits $\epsilon$-uniform mixing if and only if $X$ or $\comp{X}$ has parameters $(4m^2, 2m^2\pm m, m^2\pm m, m^2\pm m)$ where $m\ge 2$.
\end{theorem}

\section*{Acknowledgement}
This material is based upon work supported by the National Science Foundation under Grant No. 2348399.

\bibliographystyle{amsplain}
\bibliography{qw.bib}

\end{document}